\documentclass[11pt]{article}

\newcommand{\FIGDIR}{.}

\usepackage{epsf,latexsym,amssymb}
\usepackage{graphicx}
\usepackage{amsmath}  
\newcommand{\pic}[2]{\includegraphics[scale=#1]{\FIGDIR/#2}}
\newcommand{\picc}[2]{\begin{center}\pic{#1}{#2}\end{center}}
\newcommand{\edo}{\end{document}}

\newcommand{\R}{{\mathbb R}}  
\newcommand{\Z}{{\mathbb Z}}  

\newcommand{\source}{{\mathbb S}}  
\newcommand{\target}{{\mathbb T}}  

\newcommand{\order}[1]{\oplus #1}

\newcommand{\twoif}[4]{
\left\{ \begin{array}{ll}#1&#2\\#3&#4\end{array}\right.
}

\newcommand{\be}[1]{\begin{equation}\label{#1}}
\newcommand{\ee}{\end{equation}}

\newcommand{\beq}{\begin{eqnarray}}
\newcommand{\eeq}{\end{eqnarray}}
\newcommand{\beqn}{\begin{eqnarray*}}
\newcommand{\eeqn}{\end{eqnarray*}}
\newcommand{\bi}{\begin{itemize}}
\newcommand{\ei}{\end{itemize}}
\newcommand{\ben}{\begin{enumerate}}
\newcommand{\een}{\end{enumerate}}
\newcommand{\bes}[1]{\begin{subequations}\label{#1}\begin{eqnarray}}
\newcommand{\ees}[1]{\end{eqnarray}\end{subequations}}

\newcommand{\mypmatrix}[1]{\left(\begin{array}{cccccccccccc}#1\end{array}\right)}

\newcommand{\ns}{n}  
\newcommand{\Sn}{S}  
\newcommand{\nr}{m}  
\newcommand{\ncon}{q}  

\newcommand{\Rr}{R} 
\newcommand{\Rl}{\mathbf{R}} 
\newcommand{\Sl}{\mathbf{S}} 
\newcommand{\Dr}{D}  

\newcommand{\Reach}{{\cal R}}

\newcommand{\Zpn}{\Z_{\geq 0}^\ns}

\newcommand{\rates}{\rho } 

\newcommand{\ratesi}{\rates_{\indr}}

\newcommand{\mychoose}[2]{{#1 \choose #2}}
\newcommand{\Rpn}{\R_{\geq 0}^n}
\newcommand{\Rpnstrict}{\R_{>0}^n}

\newcommand{\prob}[1]{{\mathbb P}\,{\left[#1\right]}}
\newcommand{\Prob}[1]{\prob{#1}}

\newcommand{\Ex}[1]{{\mathbb E}\left[#1\right]}

\newcommand{\inds}{i}
\newcommand{\indr}{j}
\newcommand{\indl}{\pi }

\usepackage{chemarr}
\newcommand{\arrowchemnew}[1]{\xrightarrow{#1}}
\newcommand{\arrowschemnew}[2]{\xrightleftharpoons[#2]{#1}}

\newcommand{\arrowchem}[1]{\xrightarrow{#1}}

\newcommand{\Mon}{U}
\newcommand{\mon}{u}
\newcommand{\card}[1]{\#(#1)}

\newcommand{\coeflin}{\kappa _{c_j}}

\newcommand{\coefmon}{\kappa _j}
\newcommand{\coefmontilde}{{\widetilde \kappa }}

\newcommand{\Ind}{{\cal I}}  

\newcommand{\Species}{S}
\newcommand{\nocomplexes}{n_c}

\newcommand{\atopnew}[2]{\genfrac{}{}{0pt}{}{#1}{#2}}

\newcommand{\barx}{{\bar \lambda }}
\newcommand{\xstate}{\lambda }

\newcommand{\barxi}{{\bar \xi }}

\newcommand{\indsum}{\ell}

\newcommand{\bone}{\beta _1} 
\newcommand{\bonep}{p} 
\newcommand{\btwo}{\beta _2}
\newcommand{\btwon}{n} 

\newcommand{\bcon}{\beta _{\ncon}}
\newcommand{\bii}{\beta _i}
\newcommand{\amatrixZ}{\alpha } 
\newcommand{\heaviside}[1]{{\cal H}(#1)}

\newcommand{\onekaj}{\heaviside{k-a_{{\indr}}}}

\newcommand{\Exp}[1]{E[#1]}
\newcommand{\Sone}{\Exp{\Species_1}}
\newcommand{\Stwo}{\Exp{\Species_2}}
\newcommand{\Sthree}{\Exp{\Species_3}}
\newcommand{\Soneone}{\Exp{\Species_1^2}}
\newcommand{\Sthreethree}{\Exp{\Species_3^2}}
\newcommand{\Sonethree}{\Exp{\Species_1\Species_3}}
\newcommand{\Sonetwo}{\Exp{\Species_1\Species_2}}

\newcommand{\nogenes}{\beta }

\newcommand{\Speciesnew}{\Species}

\title{Examples of computation of exact moment dynamics\\
  for chemical reaction networks}
\author{Eduardo D. Sontag\\
  Rutgers University\\
{\tt eduardo.sontag@gmail.com}}

\begin{document}
\maketitle

\begin{abstract}

  \noindent
We review in a unified way results for two types of stochastic chemical reaction systems for
which moments can be effectively computed:
\ben
\item
  \emph{feedforward networks (FFN)}, treated in
  \cite{sontag_singh_june2015}, and
\item
  \emph{complex-balanced networks (CBN)}, treated in 
  \cite{mvp},
  \een
and provide several worked examples.
\end{abstract}

\section{Stochastic Kinetics}

We start by reviewing standard concepts regarding master equations for
biochemical networks, see for instance~\cite{sysbionotes}.

\subsection{Chemical Reaction Networks}

Chemical reaction networks involve interactions among a finite set of \emph{species}
\[
{\cal S} = \{\Sn_\inds , \,\inds = 1,2, \ldots \ns \}
\]
where one thinks of the $S_\inds$'s as counting the numbers of molecules
of a certain type (or individuals in an ecological model, or cells in a cell
population model):
\[
\Species_\inds(t) = k_i = \mbox{ number of units of species }
\inds\mbox{ at time }t\,.
\]
In stochastic models, one thinks of these as random variables, which interact
with each other.  The complete vector
$\Species=(\Species_1,\ldots ,\Species_\ns)'$
is called the \emph{state} of the system at time $t$, and it is
probabilistically described as a Markov stochastic process
which is indexed by time $t\geq 0$ and takes values in $\Zpn$.
Thus, $\Species(t)$ is a $\Zpn$-valued random variable, for each $t\geq 0$.
(Abusing notation, we also write $\Species(t)$ to represent an outcome of this
random variable on a realization of the process.)
We will denote
\[
p_k(t) = \prob{\Species(t) = k} 
\]
for each $k\in \Zpn$.
Then $p(t) = (p_k)_{k\in \Zpn}$ is the discrete probability density
(also called the ``probability mass function'') of $\Species(t)$.
To describe the Markov process, one needs to formally introduce chemical
reaction networks.

Mathematically, a \emph{chemical reaction network} is a finite set
\[
{\cal R} = \{\Rr_\indr, \indr = 1,2, \ldots , \nr \}
\]
of formal transformations or \emph{reactions}
\be{eq:eq:reaction}
\Rr_{\indr}: \quad \sum_{{\inds} =1}^{\ns} a_{{\inds}{\indr}} \Sn_{\inds}
\;\longrightarrow\;
 \sum_{{\inds} =1}^{\ns} b_{{\inds}{\indr}} \Sn_{\inds}
\,,\quad \indr\in \{1,2, \ldots , \nr \}
\ee
among species, together with a set of $\nr$ functions 
\[
\ratesi:\Zpn\rightarrow \R_{\geq 0}, \quad {\indr}=1,\ldots ,\nr,
\mbox{ with }\ratesi(0)=0
\]
called the \emph{propensity functions} for the respective reactions $\Rr_\indr$.
The coefficients $a_{{\inds}{\indr}}$ and $b_{{\inds}{\indr}}$ are non-negative
integers, called the  \emph{stoichiometry coefficients},
and the sums are understood informally, indicating combinations of elements.
The intuitive interpretation is 
that $\ratesi(S_1,\ldots ,S_\ns)dt$ is 
the probability that reaction $\Rr_{\indr}$ takes place, in a short interval of
length $dt$, provided that the complete state was $S=(S_1,\ldots ,S_\ns)$ at the
beginning of the interval.
In principle, the propensities can be quite arbitrary functions, but we
will focus on mass-action kinetics, for which the functions $\ratesi$ are
polynomials whose degree is the sum of the $a_{{\inds}{\indr}}$'s in the
respective reaction.  Before discussing propensities, however,
we need to introduce some more notations and terminology.

The linear combinations
$\sum_{{\inds} =1}^{\ns} a_{{\inds}{\indr}} \Sn_{\inds}$
and
$\sum_{{\inds} =1}^{\ns} b_{{\inds}{\indr}} \Sn_{\inds}$
appearing in the $\nr$ reactions are called the \emph{complexes} involved in the
reactions.
For each reaction $\Rr_\indr$, we collect the coefficients appearing on its
left-hand side and on its right-hand side into two vectors, respectively:
\[
\source(\Rr_\indr) = a_{\indr} := (a_{1{\indr}},\ldots ,a_{\ns {\indr}})'
\,,\quad
\target(\Rr_\indr) = b_{\indr} := (b_{1{\indr}},\ldots ,b_{\ns {\indr}})'
\]
(prime indicates transpose).
We call $\source,\target: {\cal R}\rightarrow {\cal C}$ the \emph{source} and \emph{target}
functions, where ${\cal C}\subseteq \Z^\ns_{\geq 0}$ is the set of all vectors
$\left\{a_{\indr},b_{\indr}, j=1\ldots \nr\right\}$.
We identify complexes with elements of ${\cal C}$.
The \emph{reactants} $S_\inds$ of the reaction $\Rr_{\indr}$ are those species
appearing with a nonzero coefficient, $a_{{\inds}{\indr}}\not= 0$ in its left-hand
side and the
\emph{products} $S_\inds$ of reaction $\Rr_{\indr}$ are those species appearing
with a nonzero coefficient $b_{{\inds}{\indr}}\not= 0$ in its right-hand side.

For every vector of non-negative integers $v=(v_1,\ldots ,v_\ns)\in \Zpn$, let us
write the sum of its entries as
$\order{v} :=v_1+\ldots +v_\ns$.
In particular, for each $\indr\in \{1,\ldots ,\nr\}$, we define the \emph{order} of
the reaction $\Rr_\indr$ as 
$
\order{a_j} :=\,
\sum_{{\inds}=1}^{\ns} 
a_{{\inds}{\indr}}
$,
which is the total number of units of all species  participating in the
reaction $\Rr_\indr$.

The $\ns\times \nr$ \emph{stoichiometry matrix} $\Gamma  = \{\gamma _{{\inds}{\indr}}\}$ is
defined as the matrix whose entries are defined as follows:
\[
\gamma _{{\inds}{\indr}} := \, b_{{\inds}{\indr}}-a_{{\inds}{\indr}} \,, \quad
{\inds}=1,\ldots ,\ns, \;{\indr}=1,\ldots ,\nr
\,.
\]
The integer $\gamma _{{\inds}{\indr}}$ counts the net change (positive or negative)
in the number of units of
species $\Sn_\inds$ each time that the reaction $\Rr_\indr$ takes place.
We will denote by $\gamma _{\indr}$ the ${\indr}$th column of $\Gamma $.  Note that,
with the notations introduced so far,
\[
\gamma _{\indr} \,=\, b_{\indr} - a_{\indr}\,, \quad
{\indr}=1,\ldots ,\nr\,.
\]
We will assume, to avoid trivial situations, that $\gamma _{\indr}\not= 0$
for all $\indr$ (that is, each reaction changes at least some species).

For example, suppose that $\ns=4$, $\nr=2$, and the reactions are
\[
R_1: \, S_1+S_2 \,\rightarrow \, S_3+S_4\,,\;\;\; R_2: \, 2S_1+S_3 \,\rightarrow \, S_2
\]
which have orders $1+1=2$ and $2+1=3$ respectively.
The set ${\cal C}$ has four elements, which list the coefficients of
the species participating in the reactions:
\[
{\cal C} = \{(1,1,0,0)', (0,0,1,1)', (2,0,1,0)', (0,1,0,0)'\}
\]
with
\[
\source(\Rr_1) = a_1 = (1,1,0,0)'
,\,
\source(\Rr_2) = a_2 = (2,0,1,0)'
,\,
\target(\Rr_1) = b_1 = (0,0,1,1)'
,\,
\target(\Rr_1) = b_2 = (0,1,0,0)'
\]
and
$\gamma _1 = (-1,-1,1,1)'$, $\gamma _2=(-2,1,-1,0)'$.
The reactants of $R_1$ are $S_1$ and $S_2$,
the reactants of $R_2$ are $S_1$ and $S_3$,
the products of $R_1$ are $S_3$ and $S_4$,
the only product of $R_2$ is $S_2$,
and the stoichiometry matrix is
\[
\Gamma  \;=\; \mypmatrix{%
  -1 & -2\cr
  -1 &  1\cr
   1 & -1\cr
   1 &  0} \,.
\]

It is sometimes convenient to write
\[
\sum_{{\inds} =1}^{\ns} a_{{\inds}{\indr}} \Sn_{\inds}
\;\arrowchemnew{\ratesi(S)}\;
\sum_{{\inds} =1}^{\ns} b_{{\inds}{\indr}} \Sn_{\inds}  
\]
to show that the propensity $\ratesi$ is associated to the reaction
$\indr$, and also to combine two reactions
$\Rr_\indr$ and $\Rr_k$ that are the reverse of each other, meaning that
their complexes are transposed:
$\source(\Rr_\indr)=\target(\Rr_k)$
and
$\source(\Rr_k)=\target(\Rr_\indr)$,
by using double arrows, like
\[
\sum_{{\inds} =1}^{\ns} a_{{\inds}{\indr}} \Sn_{\inds}
\;\arrowschemnew{\ratesi(S)}{\rates_k(S)}\;
\sum_{{\inds} =1}^{\ns} b_{{\inds}{\indr}} \Sn_{\inds} \,.
\]
When propensities are given by mass-action kinetics, as discussed below, 
one simply writes on the arrows the kinetic constants instead of the full form
of the kinetics.

\subsection{Chemical Master Equation}

A \emph{Chemical Master Equation (CME)}, which is the differential form of the
Chapman-Kolmogorov forward equation
is a system of linear differential
equations that describes the time evolution of the joint probability
distribution of the $\Species_i(t)$'s:
\be{eq:cme}
\frac{dp_k}{dt}
\;=\;
  \sum_{{\indr}=1}^\nr \ratesi(k-\gamma _{\indr})\,p_{k-\gamma _{\indr}} 
-
   \sum_{{\indr}=1}^\nr \ratesi(k)\,p_{k}
\,,\quad k\in \Zpn
\ee
where, for notational simplicity, we omitted the time argument ``$t$'' from $p$,
and the function $\ratesi$ has the property that 
$\ratesi(k-\gamma _{\indr})=0$ unless $k\geq \gamma _{\indr}$ (coordinatewise inequality).
There is one equation for each $k\in \Zpn$, so this is an infinite system of
linked equations. 
When discussing the CME, we will assume that an initial probability vector
$p(0)$ has been specified, and that there is a unique solution
of~(\ref{eq:cme}) defined for all $t\geq 0$.
(See \cite{meyntweedie1993} for existence and uniqueness results.)
A different CME results for each choice of propensity functions, a choice that
is dictated by physical chemistry considerations.  

The most commonly used propensity functions, and the ones best-justified from
elementary physical principles,
are \emph{ideal mass action kinetics} propensities, defined as follows
(see~\cite{gil00}), proportional to the number of ways in which species
can combine to form the $\indr$th source complex:
\be{eq:mass_action_propensities}
\ratesi (k) = \coefmon
\prod_{{\inds}=1}^{\ns}\mychoose{k_{\inds}}{a_{{\inds}{\indr}}}
\,\onekaj
\,\quad
{\indr} = 1,\ldots ,\nr.
\ee
where, for any scalar or vector, we denote $\heaviside{u}=1$ if $u\geq 0$
(coordinatewise) and $\heaviside{u}=0$ otherwise.
In other words, the expression can only be nonzero provided that
$k_{\inds} \geq a_{{\inds}{\indr}}$ for all $\inds=1,\ldots ,\ns$
(and thus the combinatorial coefficients are well-defined).
Observe that the expression in the right-hand side makes sense
even if $k\ngeq0$, in the following sense.
In that case, $k_{\inds}<0$ for some index $\inds$, so the factorial is not
well-defined, but on the other hand,
$k_{\inds}- a_{{\inds}{\indr}}\leq k_{\inds}<0$ implies that
$\heaviside{k-a_{{\indr}}}=0$.
So $\ratesi (k)$ can be thought of as defined by this formula for all
$k\in \Z^\ns$, even if some entries of $k$ are negative, but is zero unless
$k\geq 0$, and the combinatorial coefficients 
can be arbitrarily defined for $k\ngeq0$.
(In particular, $\ratesi(k-\gamma _{\indr})=0$ unless $k\geq \gamma _{\indr}$ in~\eqref{eq:cme}.)
The $\nr$ non-negative ``kinetic constants'' 
are arbitrary, and they
represent quantities related to the volume, shapes of the reactants, chemical
and physical information, and temperature.
The model described here assumes that temperature and volume are
constant, and that the system is well-mixed (no spatial heterogeneity).

\subsection{Derivatives of moments expressed as linear combinations of moments}

Notice that $\ratesi (k)$ can be expanded into a polynomial in which each
variable $k_i$ has an exponent less or equal to $a_{ij}$.  In other words,
\[
\ratesi (k) = \sum_{c_{\indr}\leq a_{\indr}} \coeflin \, k^{c_{\indr}}
\]
(``$\leq $'' is understood coordinatewise,
and by definition $k^{c_{\indr}} = k_1^{c_{1\indr}}\ldots k_\ns^{c_{\ns\indr}}$
and $r^0=1$ for all integers),
for suitably redefined coefficients $\coeflin$'s.

Suppose given a function $M: \Zpn\rightarrow \R$
(to be taken as a monomial when computing moments).
The expectation of the random variable $M(\Species)$ is by definition
\[
\Ex{M(\Species(t))} = \sum_{k\in \Zpn} p_k(t)\,M(k)
\,,
\]
since $p_k(t)=\Prob{\Species(t)=k}$.
Let us define, for any
$\gamma \in \Z^\ns$,
the new function
$\Delta _{\gamma }M$ given by
\[
(\Delta _{\gamma }M)(k):= M(k+\gamma )-M(k) \,.
\]
With these notations,
\be{eq:expectedM}
\frac{d}{dt} \Ex{M(\Species(t))}
\;=\;
\sum_{\indr=1}^\nr \Ex{\ratesi(\Species(t))\,\Delta _{\gamma _j}M(\Species(t))}
\ee
(see \cite{sysbionotes} for more details).
We next specialize to a monomial function:
\[
M(k) = k^{\mon} = k_1^{\mon_1} k_2^{\mon_2} \ldots  k_{\ns}^{\mon_{\ns}} 
\]
where $\mon\in \Zpn$.
There results
\[
(\Delta _{\gamma _\indr}M)(k) = \sum_{\nu \in \Ind(\mon,\indr)} \, d_{\nu }k^{\nu }
\]
for appropriate coefficients $d_{\nu }$, where
\[
\Ind(\mon,\indr)
:=
\left\{
\nu \in \Zpn
\left|
\begin{aligned}
&\nu = \mon-\mu  , \;  \mon \geq  \mu  \not=  0\\
\mu _\inds=0 &\mbox{ for each } \inds \mbox{ such that } \gamma _{{\inds}{\indr}}=0 
\end{aligned}
\right.
\right\}
\]
(inequalities ``$\geq $'' in $\Zpn$ are understood coordinatewise).
Thus, for
\eqref{eq:mass_action_propensities}:
\be{eq:momentsode}
\frac{d}{dt} \Ex{\Species(t)^{\mon}}=
\sum_{\indr=1}^\nr \, \sum_{c_{\indr}\leq a_{\indr}} \, \sum_{\nu \in \Ind(\mon,\indr)} 
d_{\nu } \coeflin \, \Ex{\Species(t)^{\nu +c_j}}.
\ee
In other words, we can recursively express the derivative of the moment of
order $\mon$ as a linear combination of other moments.  This results in an
infinite set of coupled linear ordinary differential equations, so it is natural
to ask whether, for given a particular moment or order $\mon$ of interest,
there is a finite set of moments, including the desired one, that satisfies a
finite set of differential equations.  This question can be reformulated
combinatorially, as follows.

For each multi-index $\mon\in \Zpn$, let us define $\Reach^0(\mon)=\{\mon\}$,
\[
\Reach^1(\mon) \,:= \;
\left\{\nu  + c_j \,,\; 1\leq \indr\leq \nr , c_{\indr}\leq a_{\indr} , \, \nu \in \Ind(\mon,\indr) \right\}
\]
and, more generally,
for any $\ell\geq 1$,
\[
\Reach^{\ell+1}(\mon) \,:= \; \Reach^1(\Reach^{\ell}(u))
\]
where, for any set $U$, $\Reach^\ell(\Mon) := \bigcup _{\mon\in \Mon} \Reach^\ell(\mon)$.
Finally, we set
\[
\Reach(\mon) \,:=\; \bigcup _{i=0}^{\infty } \Reach^i(\mon).
\]
Each set $\Reach^\ell(\mon)$ is finite, but the cardinality $\card{\Reach(\mon)}$ may be
infinite. 
It is finite if and only if there is some $L\geq 0$ such that
$\Reach(\mon)=\bigcup _{i=0}^{L} \Reach^i(\mon)$, or equivalently
$\Reach^{L+1}(\mon)\subseteq \bigcup _{i=0}^{L} \Reach^i(\mon)$.

Equation \eqref{eq:momentsode} says that
the derivative of the $\mon$-th moment can be expressed as a linear combination
of the moments in the set $\Reach^1(u)$.  The derivatives of these moments, in
turn, can be expressed in terms of the moments in the set $\Reach^1(u')$, for each
$u'\in \Reach^1(u)$, i.e., in terms of moments in the set $\Reach^2(u)$.
Iterating, we have the following:
``Finite reachability implies linear moment closure'' observation:

{\bf Lemma.}
Suppose that $N:=\card{\Reach(\mon)}<\infty $, and
$\Reach(u) = \{u=\mon_1,\ldots ,\mon_N\}$.
Then, writing
\[
x(t) := \left(\Ex{\Species^{\mon_1}(t)},\ldots ,\Ex{\Species^{\mon_N}(t)}\right)',
\]
there is an $A\in \R^{N\times  N}$ such that
$\dot x(t) = Ax(t)$ for all $t\geq 0$.

A classical case is that in which all reactions have order 0 or 1,
i.e. $\order{a_{\indr}}\in \{0,1\}$.
Indeed, since $\mu \not= 0$ in the definition of $\Ind(\mon,\indr)$, it follows
that $\order{a_{\indr}}\leq \order{\mu }$ for every index $\indr$.
Therefore,
$\order{(\nu +a_\indr)}=\order{\mon}+\order{a_\indr}-\order{\mu } \leq   \order{\mon}$
for all $\mon$, and the same holds for
$\nu +c_{\indr}$ if $c_{\indr}\leq a_{\indr}$.
So all elements in $\Reach(\mon)$ have degree $\leq \order{\mon}$, and thus
$\card{\Reach(\mon)}<\infty $.
A more general case is as follows.

\section{Feedforward networks}

A chemical network will be said to be of \emph{feedforward type (FFN)} if one
can partition its $\ns$ species
\[
\Sn_\inds \,, \quad \inds\in \{1,2, \ldots , \ns \}
\]
into $p$ layers
\[
\Sl_1,\ldots ,\Sl_p
\]
and there are a total of $\nr' = \nr+d$ reactions,
where $d$ of the reactions are ``pure degradation'' (or ``dilution'') reactions
\[
\Dr_\indr: \Sn_{\inds_{\indr}}\rightarrow 0\,, \quad \indr\in \{1,\ldots ,d\}
\]
and the additional $\nr$ reactions
\[
\Rr_\indr \,, \quad \indr\in \{1,2, \ldots , \nr \}
\]
can be partitioned into $p\geq 1$ layers
\[
\Rl_1,\ldots ,\Rl_p
\]
in such a manner that, in the each reaction layer $\Rr_{\pi }$
there may be any number of order-zero or order-one reactions involving species
in layer $\pi $, but every higher-order reaction at a layer $\pi >1$
must have the form:
\[
a_{{\inds_1}{\indr}} \Sn_{{\inds_1}} + \ldots  
a_{{\inds_q}{\indr}} \Sn_{{\inds_q}}
\quad\longrightarrow\quad
a_{{\inds_1}{\indr}} \Sn_{{\inds_1}} + \ldots  
a_{{\inds_q}{\indr}} \Sn_{{\inds_q}} +
b_{{\inds_{q+1}}{\indr}} \Sn_{{\inds_{q+1}}} + \ldots  
b_{{\inds_{q+q'}}{\indr}} \Sn_{{\inds_{q+q'}}}
\]
where all the species $\Sn_{{\inds_1}}, \ldots , \Sn_{{\inds_q}}$ belong to layers
having indices $<\indl$, and the species
$\Sn_{{\inds_{q+1}}},\ldots ,\Sn_{{\inds_{q+q'}}}$ are in layer $\indl$.
In other words, multimers of species in ``previous'' layers can ``catalyze''
the production of species in the given layer, but are not affected by these
reactions. 
This can be summarized by saying that for reactions at any given layer
$\indl$, the only species that appear as reactants in nonlinear reactions are
those in layers $<\indl$ and the only ones that can change are those in layer
$\indl$.

A more formal way to state the requirements is as follows.
The reactions $\Rr_\indr$ that belong to the first layer
$\Rl_1$ are all of order zero or one, i.e.\ they have
$\order{a_{\indr}}\in \{0,1\}$
(this first layer might model several independent separate chemical
subnetworks; we collect them all as one larger network), and
\be{eq:conditionslayer}
\mbox{if }\Rr_\indr\in \Rl_{\indl}:
\twoif{a_{{\inds}{\indr}}\not= 0
\mbox{ and } \order{a_\indr}>1
}%
{ \Rightarrow \Sn_\inds\in \bigcup _{1\leq s<\indl} \Sl_{\indl}}%
{\gamma _{{\inds}{\indr}}\not= 0}%
{ \Rightarrow \Sn_\inds\in \Sl_{\indl} \,.}
\ee

It was shown in~\cite{sontag_singh_june2015} that
FFN's have the finite reachability property: given any
desired moment $\mon$, there is a linear differential equation $\dot x(t) = Ax(t)$
for a suitable set of $N$ moments
\[
x(t) := \left(\Ex{\Species^{\mon_1}(t)},\ldots ,\Ex{\Species^{\mon_N}(t)}\right)',
\]
which contains the moment $u$ of interest.
Notice that steady-state moments can then be computed by solving $Ax=0$.

In practice, we simply compute 
\eqref{eq:momentsode}
starting from a desired moment, then recursively apply the same rule
to the moments appearing in the right-hand side, and so forth until no new
moments appear.  The integer $N$ at which the system closes might be very
large, but the procedure is guaranteed to stop.

We also remark that the last section of the paper~\cite{sontag_singh_june2015} explains how
certain non-feedforward networks also lead to moment closure, provided that
conservation laws ensure that variables appearing ion nonlinear reactions take
only a finite set of possible values.  We do not discuss that case here.

\subsection{Steady-states of CME}

Often, the interest is in long-time behavior, after a transient,
that is to say in the probabilistic \emph{steady state} of the system:
the joint distribution of the random variables $\Species_i=\Species_i(\infty )$ that result in the
limit as $t\rightarrow \infty $ (provided that such a limit exists in an appropriate
technical sense).
This joint distribution is a solution of the steady state CME (ssCME), the
infinite set of linear equations obtained by setting the right-hand side of
the CME to zero, that is:
\be{eq:cmess}
  \sum_{{\indr}=1}^\nr \ratesi(k-\gamma _{\indr})\,p_{k-\gamma _{\indr}} 
\;=\;
   \sum_{{\indr}=1}^\nr \ratesi(k)\,p_{k}
   \,,\quad k\in \Zpn
\ee
with the convention that $\ratesi(k-\gamma _{\indr})=0$ unless $k\geq \gamma _{\indr}$.
When substituting mass action propensities
\[
\ratesi (k) = \coefmon
\prod_{{\inds}=1}^{\ns}\mychoose{k_{\inds}}{a_{{\inds}{\indr}}}
\,\onekaj
\]
the steady-state equation~\eqref{eq:cmess} becomes:
\be{eq:cmessmassaction}
\sum_{{\indr}=1}^\nr
\,\coefmon\,
\prod_{{\inds}=1}^{\ns}\mychoose{k_{\inds}-\gamma _{\inds\indr}}{a_{{\inds}{\indr}}}
\,\heaviside{k-b_{\indr}}
\,p_{k-\gamma _{\indr}} 
\;=\;
\sum_{{\indr}=1}^\nr
\,\coefmon\,
\prod_{{\inds}=1}^{\ns}\mychoose{k_{\inds}}{a_{{\inds}{\indr}}}
\,\onekaj
\,p_{k}
   \,,\quad k\in \Zpn\,.
\ee
Equivalently,
\be{eq:cmessmassaction1}
\sum_{{\indr}=1}^\nr
\,\coefmontilde_\indr\,
\prod_{{\inds}=1}^{\ns}
\frac{\left(k_{\inds}-\gamma _{\inds\indr}\right)!}{\left(k_{\inds}-b_{{\inds}{\indr}}\right)!}
\,\heaviside{k-b_{\indr}}
\,p_{k-\gamma _{\indr}} 
\;=\;
\sum_{\indr=1}^\nr
\,\coefmontilde_\indr\,
\prod_{{\inds}=1}^{\ns}
\frac{k_{\inds}!}{\left(k_{\inds}-a_{{\inds}{\indr}}\right)!}
\,\onekaj
\,p_{k}
   \,,\quad k\in \Zpn
\ee
when introducing the new constants
$
\coefmontilde_\indr :=
{\coefmon}/
{\prod_{{\inds}=1}^{\ns}\left(a_{{\inds}{\indr}}!\right)}
$.
If, for convenience, we write:
\[
{\xstate}^k \,:=\; \xstate_1^{k_1}\ldots \xstate_\ns^{k_\ns}
\quad\mbox{and}\quad
k! \,:= \;k_1! \ldots  k_\ns!
\]
for each nonnegative integer vector $k=(k_1,\ldots ,k_\ns)\in \Zpn$
and positive vector $\xstate = (\xstate_1,\ldots ,\xstate_\ns)\in \Rpnstrict$,
then~\eqref{eq:cmessmassactionrev} is re-written as
\be{eq:cmessmassactionrev}
\sum_{{\indr}=1}^\nr
\,\coefmontilde_\indr\,
\frac{\left(k-\gamma _{\indr}\right)!}{\left(k-b_{\indr}\right)!}
\,\heaviside{k-b_{\indr}}
\,p_{k-\gamma _{\indr}} 
\;\;=\;\;
\sum_{{\indr}=1}^\nr
\,\coefmontilde_\indr\,
\frac{k!}{\left(k-a_{\indr}\right)!}
\,\onekaj
\,p_{k}
   \,,\quad k\in \Zpn
\ee
Since~\eqref{eq:cmessmassactionrev} is a linear equation on the
$\left\{p_k, k\in \Zpn \right\}$, any rescaling of a given set of 
$p_k$'s will satisfy the same equation; for probability densities, one
normalizes to a unit sum.

If there are conservation laws satisfied by the system
then steady state solutions will not be unique,
and the equation $Ax=0$ must be supplemented by a set of linear constraints that
uniquely specify the solution.
This is obvious even for simple linear reactions.
For example, suppose we consider one reversible reaction
\[
\Species_1\;\arrowschemnew{\kappa _1}{\kappa _2}\; \Species_2
\]
(propensities are mass-action, $\rates_i(\Species_1,\Species_2)=\kappa _i\Species_i$).
The first moments (means) satisfy
$\dot x_1=\kappa _2x_2-\kappa _1x_1$ and $\dot x_2=\kappa _1x_1-\kappa _2x_2$.
Any vector $(\barxi_1,\barxi_2)$ with $\kappa _1\barxi_1=\kappa _2\barxi_2$ is a steady
state of these equations. 
However, the sum of the numbers of molecules $\Species_1$ and $\Species_2$
is conserved in the reactions.
Given a particular total number, $\beta $,
the differential equations can be reduced to just one equation, say for $x_1$:
$\dot x_1 = \kappa _2(\beta -x_1)-\kappa _1x_1 = -(\kappa _1+\kappa _2)x_1 + \kappa _2\beta $, which has the
affine form $\dot x=Ax+b$.
At steady state, we have the unique solution
$\barxi_1 = \beta \kappa _2/(\kappa _1+\kappa _2)$,
$\barxi_2 = \beta \kappa _1/(\kappa _1+\kappa _2)$
obtained by imposing the constraint $\barxi_1+\barxi_2=\beta $.
It can easily be proved (see e.g.~\cite{mvp}) that at
steady state, $\Species_1$ is a binomial random variable $B(\beta ,p)$
with $p=\frac{1}{1+\mu }$, where $\mu =\kappa _1/\kappa _2$.
We later discuss further the effect of conservation laws.

\subsection{A worked example}

For networks with only zero and first order reactions, which are trivially in
our feedforward class, it is well-known that one may compute all moments in
closed form.  We work out a simple example here.
Again we start with a reversible reaction
\[
\Species_1\;\arrowschemnew{\kappa }{\delta }\; \Species_2
\]
and mass-action propensities.
We think in this context of $\Species_1$ as the active form of a certain gene
and $\Species_2$ as the inactive form of the same gene.
Transcription and translation are summarized, for simplicity, as one
reaction
\[
\Species_1\;\arrowchemnew{\rho }\; \Species_1 + \Species_3
\]
and degradation or dilution of the gene product $\Species_3$ is a linear
reaction
\[
\Species_3\;\arrowchemnew{\eta }\; \emptyset \,.
\]
The stoichiometry matrix is
\[
  \Gamma  = \mypmatrix{ -1 & 1 & 0 & 0\cr 
	 1 & -1 & 0 & 0 \cr 
         0 &  0 & 1 & -1}\,.
\]
Suppose that we are interested in the mean and variance of $\Species_3$ subject
to the conservation law $\Species_1+\Species_2=\beta $, for some fixed positive
integer $\beta $.
A linear differential equation for the vector of second order moments:
\[
{\cal M} \; = \;(\Sone, \Soneone, \Sonethree, \Sthree, \Sthreethree)'
\]
is $\dot {\cal M} = A{\cal M} + b$, where
\[
A \;=\;
\left(\begin{array}{ccccc}  - \delta  - \kappa  & 0 & 0 & 0 & 0\\ \kappa  - \delta  + 2\, \delta \, \nogenes &  - 2\, \delta  - 2\, \kappa  & 0 & 0 & 0\\ 0 & \rho  &  - \delta  - \eta  - \kappa  & \delta \, \nogenes & 0\\ \rho  & 0 & 0 & - \eta  & 0\\ \rho  & 0 & 2\, \rho  & \eta  & - 2\, \eta  \end{array}\right)
\;,\quad\quad
b \;=\; \left(\begin{array}{c}  \delta \, \nogenes\\  \delta \, \nogenes\\ 0\\ 0\\ 0 \end{array}\right)
\]

Solving $A\overline{{\cal M}}+b=0$ to obtain steady state moments, we have:
\[
\overline{{\cal M}}
\;=\;
\left(\begin{array}{c} \frac{\delta \, \nogenes}{\delta  + \kappa }\\ \frac{\delta \, \nogenes\, \left(\kappa  + \delta \, \nogenes\right)}{{\left(\delta  + \kappa \right)}^2}\\ \frac{\delta \, \nogenes\, \rho \, \left(\eta \, \kappa  + {\delta }^2\, \nogenes + \delta \, \eta \, \nogenes + \delta \, \kappa \, \nogenes\right)}{\eta \, {\left(\delta  + \kappa \right)}^2\, \left(\delta  + \eta  + \kappa \right)}\\ \frac{\delta \, \nogenes\, \rho }{\eta \, \left(\delta  + \kappa \right)}\\ \frac{\delta \, \nogenes\, \rho \, \left({\delta }^2\, \eta  + \nogenes\, \rho \, {\delta }^2 + \delta \, {\eta }^2 + 2\, \delta \, \eta \, \kappa  + \nogenes\, \rho \, \delta \, \eta  + \nogenes\, \rho \, \delta \, \kappa  + {\eta }^2\, \kappa  + \eta \, {\kappa }^2 + \rho \, \eta \, \kappa \right)}{{\eta }^2\, {\left(\delta  + \kappa \right)}^2\, \left(\delta  + \eta  + \kappa \right)} \end{array}\right)\,.
\]

\subsection{A simple nonlinear example}

We consider a feedforward system in which there are three species;
$\Species_1$ catalyzes production $\Species_2$, and $\Species_1$ and
$\Species_2$ are both needed to produce $\Species_3$:
\[
0 \arrowchemnew{\kappa _1} S_1 \arrowchemnew{\delta _1} 0
\quad
S_1 \arrowchemnew{\kappa _2} S_1+S_2
\quad
S_2 \arrowchemnew{\delta _2} 0
\quad
S_1+S_2 \arrowchemnew{\kappa _3} S_1+S_2+S_3
\quad
S_3 \arrowchemnew{\delta _3} 0
\]
Computing $\Sthree$, the mean of $\Species_3$, requires a minimal differential
equation of order 5, for the moments
\[
{\cal M} \; = \; (\Sthree, \Sonetwo, \Stwo, \Soneone, \Sone)'
\]
and has form
$\dot {\cal M} = A{\cal M} + b$, where
\[
A = \left(\begin{array}{ccccc} - \delta _2 & \kappa _3 & 0 & 0 & 0\\ 0 &  - \delta _1 - \delta _2 & \kappa _1 & \kappa _2 & 0\\ 0 & 0 & - \delta _2 & 0 & \kappa _2\\ 0 & 0 & 0 & - 2\, \delta _1 & 2\, \kappa _1 + \delta _1\\ 0 & 0 & 0 & 0 & - \delta _1 \end{array}\right)
\quad\mbox{and}\quad
b = \left(\begin{array}{c} 0\\ 0\\ 0\\ \kappa _1\\ \kappa _1 \end{array}\right) \,.
\]

\section{Poisson-like solutions and complex-balanced networks}

We observe that for any given positive vector $\barx\in \Rpnstrict$,
the set of numbers
\be{eq:equilexpo}
\Pi  \,= \,\left\{p_k = \frac{\barx^k}{k!} \;,\;\; k\in \Zpn \right\}
\ee
satisfies the ssCME equations~\eqref{eq:cmessmassactionrev} if and only if
\be{eq:cmessmassactionlambda}
\sum_{{\indr}=1}^\nr
\,\coefmontilde_\indr\,
\frac{\barx^{k-\gamma _{\indr}}}%
     {\left(k-b_{\indr}\right)!}
\,\heaviside{k-b_{\indr}}
\;\;=\;\;
\sum_{{\indr}=1}^\nr
\,\coefmontilde_\indr\,
\frac{\barx^k}{\left(k-a_{\indr}\right)!}
\,\onekaj
   \,,\quad k\in \Zpn
   \,,
\ee
Re-writing this as:
\be{eq:cmessmassactionlambdarearranged}
\sum_{c\in {\cal C}} \sum_{\{j|b_{\indr}=c\}}
\,\coefmontilde_\indr\,
\frac{\barx^{k-\gamma _{\indr}}}%
     {\left(k-b_{\indr}\right)!}
\,\heaviside{k-b_{\indr}}
     \;\;=\;\;
\sum_{c\in {\cal C}} \sum_{\{j|a_{\indr}=c\}}
\,\coefmontilde_\indr\,
\frac{\barx^k}{\left(k-c\right)!}
\,\onekaj
   \,,\quad k\in \Zpn
   \,,
\ee
we see that a sufficient condition for~\eqref{eq:equilexpo} to be a solution
is that \emph{for each individual} complex $c\in {\cal C}$:
\[
\sum_{\{j|b_{\indr}=c\}}
\,\coefmontilde_\indr\,
\frac{\barx^{k-\gamma _{\indr}}}%
     {\left(k-c\right)!}
\,\heaviside{k-b_{\indr}}
\;\;=\;\;
\sum_{\{j|a_{\indr}=c\}}
\,\coefmontilde_\indr\,
\frac{\barx^k}{\left(k-c\right)!}
\,\onekaj
\,,\quad k\in \Zpn
   \,,
\]
or, equivalently,
\[
\frac{\heaviside{k-c}}%
     {\left(k-c\right)!}
\sum_{\{j|b_{\indr}=c\}}
\,\coefmontilde_\indr\,
     \barx^{k-\gamma _{\indr}}
\;\;=\;\;
\frac{\heaviside{k-c}}{\left(k-c\right)!}
\sum_{\{j|a_{\indr}=c\}}
\,\coefmontilde_\indr\,
\barx^k
\,,\quad k\in \Zpn
   \,.
   \]
A sufficient condition for this to hold is that
\be{eq:eachcomplex}
\sum_{\{j|b_{\indr}=c\}}
\,\coefmontilde_\indr\,
\barx^{a_\indr}
\;\;=\;\;
\sum_{\{j|a_{\indr}=c\}}
\,\coefmontilde_\indr\,
\barx^{a_{\indr}}
   \,,\quad k\in \Zpn
\ee
holds for all complexes
(and conversely, this last condition is necessary for all complexes for which
$k\geq c$). 
Note that one can equally well write ``${\barx}^{c}$'' and bring this term
outside of the sum, in the right-hand side.

\begin{figure}[ht]
  \picc{0.3}{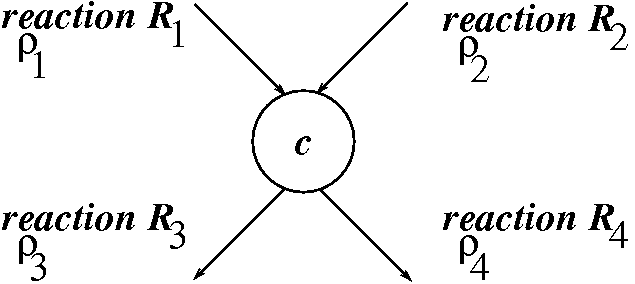}
  \caption{Complex balancing: outflows and inflows must balance at each
    complex $c$.
    The left-hand side of \protect{\eqref{eq:eachcomplex}} is
$\coefmontilde_3\barx^{a_3}+\coefmontilde_4\barx^{a_4}$
    and the right-hand side is
$\coefmontilde_1\barx^{a_1}+\coefmontilde_2\barx^{a_2}$
      }
  \label{fig:complexbal}
\end{figure}

When property~\eqref{eq:eachcomplex} holds for every complex, one says that
$\barx$ is a \emph{complex balanced steady state} of the associated
\emph{deterministic} chemical reaction network.
(That is, the system of differential equations
$\dot x = \Gamma  Q(x)$, where $Q(x)$ is a column vector of size $\nr$ whose
$\indr$th entry is $\ratesi(x)$ and $x(t)\in \Rpn$ for all $t$.)
Complex balancing means that, for each complex, outflows and inflows balance
out. 
This is a Kirschoff current law (in-flux = out-flux, at each node).
See Figure~\ref{fig:complexbal}.

Foundational results in deterministic chemical network theory were obtained by 
Horn, Jackson, and Feinberg (see~\cite{feinberg1,feinberg2}).
One of the key theorems is that a
sufficient condition for the existence of a complex balanced steady state
is that the network be \emph{weakly reversible} and have \emph{deficiency zero}.
The deficiency is computed as $\nocomplexes-\ell-r$, where $\nocomplexes$ is
the number of complexes, $r$ is the rank of the matrix $\Gamma $, and $\ell$ is the
number of ``linkage classes'' (connected components of the reaction graph).
Weak reversibility means that each connected component of the reaction graph
must be strongly connected.
One of the most interesting features of this theorem is that no assumptions
need to be made about the kinetic constants.  (Of course, the choice of
the vector $\barx$ will depend on kinetic constants.)
We refer the reader to the citations for details on deficiency theory, as well
as, of interest in the present context, several examples discussed
in~\cite{mvp}. 
The theorems for weakly reversible deficiency zero networks are actually
far stronger, and they show that every possible steady state of the
corresponding deterministic network is complex balanced, and also that they
are asymptotically stable relative to stoichiometry classes.
The connection with ssCME solutions was a beautiful observation made
in~\cite{anderson10}, but can be traced to the
``nonlinear traffic equations'' from queuing theory, described in Kelly's
textbook~\cite{kelly}, Chapter 8 (see also~\cite{mairesse09} for a discussion),

The elements of $\Pi $ given by formula~\eqref{eq:equilexpo} add up to:
\[
\sum_{k\in \Zpn} p_k \;=\;
\sum_{k_1=0}^\infty \ldots \sum_{k_\ns=0}^\infty  \,
\frac{\barx_1^k}{k_1!}\ldots \frac{\barx_\ns^k}{k_\ns!} \;=\;
Z := e^{\barx_1}\ldots e^{\barx_\ns}
\]
Thus, normalizing by the total, $\{p_k/Z, k\in \Zpn\}$ is a probability
distribution.
However, because of stoichiometric constraints, solutions are typically not
unique, and general solutions appear as convex combinations of
solutions corresponding to invariant subsets of states.
A solution with only a finite number of nonzero $p_k$'s will then have a
different normalization factor $Z$.

\subsection{Conservation laws, complex-balanced case}

When steady states do not form an irreducible Markov chain, the solutions of
the form~\eqref{eq:equilexpo} are not the only solutions in the complex
balanced case.  Restrictions to each component of the Markov chain are also
solutions, as are convex combinations of such restrictions.
To formalize this idea, suppose that there is some subset ${\cal Z}_0\subseteq \Z^n$ with
the following stoichiometric invariance property:
\be{eq:stoichinvariant}
k \in  {\cal Z}_0 \;\Rightarrow \; k\pm \gamma _{\indr} \in {\cal Z}_0  \quad \forall \, \indr = 1, \ldots , \nr \,.
\ee
(The same property is then true for the complement of ${\cal Z}_0$.)
Consider the set ${\cal Z}:={\cal Z}_0\bigcap \Zpn$.
For each $k\in {\cal Z}$, the left-hand side term
in equation~\eqref{eq:cmessmassactionlambda} either involves an index
$k-\gamma _\indr>0$, and hence also in ${\cal Z}$, or
it is zero (because $k-b_\indr\geq 0$ implies $k-\gamma _\indr\geq 0$) and so
it does not matter that $k-\gamma _\indr\not\in {\cal Z}$.
Thus,
\be{eq:equilexpo_reduced}
p_k = \twoif%
{\displaystyle
  \frac{\barx^k}{k!}}%
  {\mbox{if } k \in  {\cal Z}}%
  {0}%
  {\mbox{if } k \in  \Zpn \setminus {\cal Z}}
\ee
is also a solution, in the complex balanced case
(observe that, for indices in $\Zpn \setminus {\cal Z}$,
equation~\eqref{eq:cmessmassactionlambda} is trivially satisfied, since both
sides vanish).
So we need to divide by the sum $Z$ of the elements
in~\eqref{eq:equilexpo_reduced} in order to normalize to a probability
distribution. 
The restriction to ${\cal Z}$ will the unique steady state distribution provided that the
restricted Markov chain has appropriate irreducibility properties.

In particular, suppose that ${\cal A}=(\alpha _{ij})\in \R^{\nr\times \ns}$ is a matrix whose nullspace includes ${\cal C}$
(for example, ${\cal A}$ could be the orthogonal complement of the ``stoichiometric
subspace'' spanned by ${\cal C}$), and pick any vector
$\beta =(\beta _1,\ldots ,\beta _\ncon)'\in \R^\ncon$. 
Then ${\cal Z}_0=\{k|\, {\cal A} k=\beta \}$ satisfies~\eqref{eq:stoichinvariant}.
In this case, the sum of the elements in~\eqref{eq:equilexpo_reduced} is:
\[
Z(\bone, \ldots , \bcon)=
\sum_{
\atopnew{k_1, \ldots , k_\ns \geq 0}%
{{\cal A}k=\beta }
}
\frac{\lambda _1^{k_1}}{k_1!}\frac{\lambda _2^{k_2}}{k_2!} \ldots   
\frac{\lambda _\ns^{k_\ns}}{k_\ns!}
\]
(value is zero if the sum is empty).
The normalized form of~\eqref{eq:equilexpo_reduced} has $p_k=0$ for
$k \in  \Zpn \setminus {\cal Z}$, and
\be{pks_normalized}
p_k \;=\; \frac{1}{Z(\bone, \ldots , \bcon)} \; \frac{\lambda _1^{k_1}}{k_1!}\frac{\lambda _2^{k_2}}{k_2!}
\ldots   \frac{\lambda _\ns^{k_\ns}}{k_\ns!}
\ee
for $k \in  {\cal Z}$.
A probabilistic interpretation of these numbers is as follows.

Suppose that we have $n$ independent Poisson random variables,
$\Speciesnew_\inds$, $\inds=1, \ldots , \ns$, with parameters
$\lambda _\inds$ respectively, so
\begin{equation}
\label{eq:P1}
\Prob{\Speciesnew_1=k_1, \Speciesnew_2=k_2, \ldots  , \Speciesnew_\ns=k_\ns} \;=\;
e^{-(\lambda _1 + \ldots  + \lambda _\ns)}
\, \frac{\lambda _1^{k_1}}{k_1!}\frac{\lambda _2^{k_2}}{k_2!}
\ldots   \frac{\lambda _\ns^{k_\ns}}{k_\ns!} 
\end{equation}
for $k\geq 0$ (and zero otherwise).
Let us introduce the following new random variables:
\[
Y_\indr \,:= \;
\sum_{\inds=1}^{n} \amatrixZ_{\indr \inds} \Speciesnew_\inds \quad , \quad \indr=1 , \ldots  , \ncon\,.
\]
Observe that
\beqn
\Prob{Y_1=\bone, \ldots , Y_m=\bcon}
&=&
\sum_{%
\atopnew{k_1, \ldots , k_\ns \geq 0}%
{\amatrixZ_{11} k_1 + \ldots  + \amatrixZ_{1n} k_\ns=\bone , \; \ldots  , \;
\amatrixZ_{\ncon 1} k_1 + \ldots  + \amatrixZ_{\ncon n} k_\ns=\bcon}}\;
\Prob{\Speciesnew_1=k_1, \Speciesnew_2=k_2, \ldots  , \Speciesnew_\ns=k_\ns}
\\
&=& 
e^{-(\lambda _1 + \ldots  + \lambda _\ns)}Z(\bone, \ldots , \bcon) \,.
\eeqn
Therefore, for each $k\in {\cal Z}$, $p_k$ in~\eqref{pks_normalized} equals the
following conditional probability:
\[
\frac{\Prob{\Speciesnew_1=k_1, \Speciesnew_2=k_2, \ldots  , \Speciesnew_\ns=k_\ns}}%
  {\Prob{Y_1=\bone, \ldots , Y_\ncon=\bcon}}
  \;=\;
  \Prob{\Speciesnew_1=k_1, \Speciesnew_2=k_2, \ldots  , \Speciesnew_\ns=k_\ns
        \, \bigl \vert \, Y_1=\beta _1, \ldots  , Y_\ncon=\beta _\ncon}\,.
  \]
If our interest is in computing this conditional probability, the main effort
goes into computing $Z(\bone, \ldots , \bcon)$.
The main contribution of the paper~\cite{mvp} was to provide effective
algorithms for the computation of
$Z(\bone, \ldots  , \bcon)$ recursively on the $\bii$'s.
A package for that purpose, called MVPoisson, was included with that paper.

Conditional moments
\[
E[\Species_j^r \, \bigl \vert \, Y_1=\bone, \ldots  , Y_m=\bcon]
\,,\; r \geq 1\,,
\]
including the conditional expectation (when $r=1$),
as well as centered moments such as the conditional variance,
can be computed once that these conditional probabilities are known. 
It is convenient for that purpose to first compute the factorial moments.
Recall that the $r$th factorial moment $E[W^{(r)}]$
of a random variable $W$ is defined as the expectation of $W!/(W-r)!$.
For example, when
$r=1$, $E[W^{(r)}] = E[W]$, and for
$r=2$, $E[W^{(r)}] = E[W^2] - E[W]$,
and thus the mean and variance of $W$ can be obtained from these.
We denote the conditional factorial moment of $\Species_\inds$ given $Y=\beta $,
as $E[\Species_j^{(r)}\, \bigl \vert \, Y]$.
It is not difficult to see (Theorem 2 in~\cite{mvp}) that:
\[
E[\Species_j^{(r)}\, \bigl \vert \, Y]
\;\;=\;\;
\lambda _j^r \, \cdot \,
\frac{Z(\bone-r\amatrixZ_{1j},\btwo- r\amatrixZ_{2j}, \ldots ,
  \bcon-r\amatrixZ_{\ncon j}) }%
{Z(\bone, \ldots , \bcon)}
\]
when all $\bii-r\amatrixZ_{ij}\geq 0$ and zero otherwise.
The paper~\cite{mvp} also discusses how mixed moments such as covariances can be
computed.

For example, from $r=1$ we obtain a formula for the conditional mean: 
\be{eq:condmean}
E[\Species_j\, \bigl \vert \, Y]
\;\;=\;\;
\lambda _j \, \cdot \,
\frac{Z(\bone-\amatrixZ_{1j},\btwo- \amatrixZ_{2j}, \ldots ,
  \bcon-\amatrixZ_{\ncon j}) }%
{Z(\bone, \ldots , \bcon)}
\ee
when all $\bii\geq \amatrixZ_{ij}$, and zero otherwise,
and from $r=2$ we obtain a formula for the conditional second moment:
\[
E[\Species_j^2\, \bigl \vert \, Y]
\;\;=\;\;
\lambda _j^2 \, \cdot \,
\frac{Z(\bone-2\amatrixZ_{1j},\btwo-2\amatrixZ_{2j}, \ldots ,
  \bcon-2\amatrixZ_{\ncon j}) }%
     {Z(\bone, \ldots , \bcon)}
\;\;+\;\; E[\Species_j\, \bigl \vert \, Y]
\]
when all $\bii\geq 2\amatrixZ_{ij}$, and zero otherwise.
We next work out a concrete example.

\subsection{A worked example}

Suppose that two molecules of species $\Species_1$ and $\Species_2$ can
reversibly combine through a bimolecular reaction to produce a molecule of
species $\Species_3$
\[
S_1+S_2\;\arrowschemnew{\kappa _1}{\kappa _2}\; S_3 \,.
\]
Since the deficiency of this network is $\nocomplexes-\ell-r=2-1-1=0$ and it is
reversible and hence weakly reversible as well, we know that there is a
complex-balanced equilibrium (and every equilibrium is complex balanced).
We may pick, for example, $\barx=(1,1,K)$, where $K := \kappa _1 / \kappa _2$.
The count of $\Species_1$ molecules goes down by one every time that a
reaction takes place, at which time the count of $\Species_3$ molecules goes
up by one.  Thus, the sum of the number of $\Species_1$ molecules plus the
number of $\Species_3$ molecules remains constant in time, equal to their
starting value, which we denote as $\bonep$.
Similarly, the sum of the number of $\Species_2$ molecules plus the
number of $\Species_3$ molecules remains constant, equal to some number $\btwon$.
(In the general notations, we have $a_{11}=a_{13}=1$, $a_{22}=a_{23}=1$,
$a_{12}=a_{21}=0$, $\beta _1=p$, $\beta _2=n$.)
In the steady state limit as $t\rightarrow \infty $, these constraints persist.
In other words, all $p_k$ should vanish except those corresponding to
vectors $k=(k_1,k_2,k_3)$ such that $k_1+k_3=\bonep$ and $k_2+k_3=\btwon$.
The set consisting of all such vectors is invariant, so
\[
p_k = \twoif%
{\displaystyle
  \frac{\barx_1^{k_1}}{k_1!}
  \frac{\barx_2^{k_2}}{k_2!}
  \frac{\barx_3^{k_3}}{k_3!}%
    }%
  {\mbox{if } k_1+k_3=\bonep \mbox{ and } k_2+k_3=\btwon}%
  {0}%
  {\mbox{otherwise}}
\]
is a solution of the ssCME.  In order to obtain a probability density, we
must normalize by the sum $Z(\bonep,\btwon)$ of these $p_k$'s.
Because of the two constraints, the sum can be expressed in terms of just
one of the indices, let us say $k_1$.
Observe that, since $k+k_3=p$ and $k_3\geq 0$, necessarily $k\leq p$.
Since $k_2=n-k_3 = n+k-p$ must be non-negative, we also have the constraint
$k\geq \max\{0,p-n\}$.
So the only nonzero terms are for $k\in \{\max\{0,p-n\},\ldots ,p\}$.
With $k_3=p-k$, $k_2=n-k_3 = n+k-p$, we have:
\be{eq:Zexample_first}
Z(p,n) \;=\;
\sum_{\indsum=\max\{0,p-n\}}^{p} \, \frac{K^{p-\indsum}} { \indsum!\, (n+\indsum-p)!\, (p-\indsum)!}
\;=\;
\sum_{\indsum=0}^{\min\{p,n\}} \, \frac{K^\indsum} { (p-\indsum)!\, (n-\indsum)!\, \indsum! }
\ee
The second form if the summation makes it obvious that $Z(p,n)=Z(n,p)$.

When $n\geq p$, we can also write
\be{eq:Zexample}
Z(p,n) \;=\;
\frac{1}{n! p!}\,
\sum_{\indsum=0}^{p} \,\, \frac{n!}{(n-p+\indsum)!} {p \choose \indsum} K^{p-\indsum}
\ee
which shows the expression as a rational function in which the numerator
is a polynomial of degree $p$ on $n$.
This was derived assuming that $n\geq p$, and the factorials in the denominator do
not make sense otherwise.
However, let us think of each term $\frac{n!}{(n-p+\indsum)!}$ as the product
$n(n-1)\ldots (n-p+\indsum+1)$, which may include zero as well as negative numbers.
With this understanding,
the formula in~\eqref{eq:Zexample} makes sense even when $n<p$.
Observe that such a term vanishes for any index $\indsum<p-n$.
Thus, for $n<p$, \eqref{eq:Zexample} reduces to:
\[
\frac{1}{p!}\,
\sum_{\indsum=p-n}^{p} \,\, \frac{1}{(n-p+\indsum)!} {p \choose \indsum} K^{p-\indsum}
\]
or equivalently, with a change of indices $~\indsum=p-\indsum$ and then using
${p \choose p-\indsum} = {p \choose \indsum}$:
\[
\frac{1}{p!}\,
\sum_{\indsum=0}^{n} \,\, \frac{1}{(n-\indsum)!} {p \choose p-\indsum} K^\indsum
\;=\;
\frac{1}{p!}\,
\sum_{\indsum=0}^{n} \,\, \frac{1}{(n-\indsum)!} {p \choose \indsum} K^\indsum
\;=\;
\sum_{\indsum=0}^{n} \,\, \frac{K^\indsum}{(n-\indsum)! (p-\indsum)! \indsum!}
\,.
\]
In this last form, we have the same expression as the last one in
\eqref{eq:Zexample_first}.
In conclusion, provided that we interpret the quotient of combinatorial numbers
in~\eqref{eq:Zexample} as a product that may be zero,
formula~\eqref{eq:Zexample} is valid for all $n$ and $p$, not just for
$n\geq p$.  
In particular, we have;
\[
Z(0,n)= \frac{1} {n!}\,,
\]
\[
Z(1,n) \;=\;
\frac{Kn+1} {n!}\,,
\]
\[
Z(2,n) \;=\;
\frac{{K}^{2}{n}^{2}+ \left( -{K}^{2}+2\,K \right) n + 1} {2n!}\,,
\]
\[
Z(3,n) \;=\;
\frac{{K}^{3}{n}^{3}+ \left( -3\,{K}^{3}+3\,{K}^{2} \right) {n}^{2}+
 \left( 2\,{K}^{3}-3\,{K}^{2}+3\,K \right) n + 1} {3!n!}
\]
and so forth.
In terms of the Gauss's hypergeometric function $_2$F$_0$, we can also write:
\[
Z(p,n) \;=\;  {\frac {{\mbox{$_2$F$_0$}(-n,-p;\,\ ;\,K)}}{p!\,n!}}
\]
The recursion on $n$ obtained by using the package MVPoisson from~\cite{mvp}
is as follows (by symmetry, a recursion on $p$ can be found by exchanging $n$
and $p$): 
\[
Z ( p,n+2 ) \;=\;
{\frac{K}{n+2}} \,Z(p,n) \,+\,{\frac{-Kn+Kp-K+1}{n+2}} \,Z(p,n+1)\,.
\]

Now~\eqref{eq:condmean} gives the conditional mean of the first species,
$\Species_1$ ($j=1$ for this index, $r=1$ for the first moment, and
$\lambda _1^1 =1^1 = 1$) as: 
\[
\varphi(p,n) \, := \;
E[\Species_1\, \bigl \vert \, \Species_1+\Species_3=p, \Species_2+\Species_3=n] \;=\;
\frac{Z(p-1,n)}%
{Z(p,n)} \,.
\]
for $p\geq 1$, $n\geq 0$, and zero otherwise.  For example,
\[
\varphi(1,n) \;=\; \frac{1}{Kn+1}
\]
\[
\varphi(2,n) \;=\;
\frac{2 (Kn+1)}{{K}^{2}{n}^{2}+ \left( -{K}^{2}+2\,K \right) n + 1}\,.
\]

\subsection{Another worked example}

Suppose that molecules of species $\Species_1$ can be randomly created and
degraded, and they can also reversibly combine with molecules of $\Species_2$
through a bimolecular reaction to produce molecules of
species $\Species_3$:
\[
\emptyset \;\arrowchem{\kappa _1} \;\Species_1 \;\arrowchem{\kappa _2} \;\emptyset\,,\quad
\Species_1+\Species_2\;\arrowschemnew{\kappa _3}{\kappa _4}\; \Species_3 \,.
\]
There are $\nocomplexes=4$ complexes: $\emptyset$, $\Species_1$,
$\Species_1+\Species_2$, and $\Species_3$, and $\ell=2$ linkage classes.
The stoichiometry matrix
\[
\Gamma  \;=\;
\mypmatrix{%
   1 &-1 & -1 & 1\cr
   0 & 0 & -1 & 1\cr
   0 & 0 &  1 & -1}
  \]
  has rank $r=2$, so the deficiency of this weakly reversible network is
$\nocomplexes-\ell-r=4-2-2=0$.
Thus there is a complex-balanced equilibrium (and every equilibrium is complex
balanced). 
We may pick, for example, $\barx=(\lambda ,1,\mu )$, where $\lambda  := \frac{\kappa _1}{\kappa _2}$
and $\mu := \frac{\kappa _1\kappa _3}{\kappa _2\kappa _4}$.
Notice that there is only one nontrivial conserved quantity,
$\Species_2+\Species_3=n$, since $\Species_1$ is not conserved.
We have:
\[
Z(n)
\;=\;
\sum_{
\atopnew{k_1, k_2, k_3 \geq 0}%
{k_2+k_3=n}
}
\frac{\lambda _1^{k_1}}{k_1!}\frac{\lambda _2^{k_2}}{k_2!} \frac{\lambda _3^{k_3}}{k_3!}
\;=\;
\sum_{k_1=0}^{\infty }\frac{\lambda ^{k_1}}{k_1!}\,
\sum_{k_2=0}^n\frac{\mu ^{n-k_2}}{k_2!(n-k_2)!}
\;=\;
\frac{e^{\lambda }}{n!} (1+\mu )^n\,.
\]
The normalized probability~\eqref{eq:equilexpo_reduced}, for
$k=(k_1,k_2,k_3)\geq 0$ with $k_2+k_3=n$, is:
\[
p_k
\;=\;
\frac{1}{Z(n)} \; \frac{\lambda ^{k_1}}{k_1!}\frac{1}{k_2!}
\frac{\mu ^{k_3}}{k_3!}
\;=\;
\frac{n!}{e^{\lambda }(1+\mu )^n}\frac{\lambda ^{k_1}\mu ^{k_3}}{k_1!k_2!k_3!}
\]
and as discussed earlier, this is the conditional probability
\[
\Prob{\Species_1=k_1, \Speciesnew_2=k_2, \Speciesnew_3=k_3 \, \bigl \vert \,
  \Species_2+\Species_3=n}\,.
\]
Using this expression, we may compute, for example, the conditional marginal
distribution of $\Species_2$: 
\[
\Prob{\Speciesnew_2=r\, \bigl \vert \,  \Species_2+\Species_3=n}
\;=\;
\sum_{k_1=0}^{\infty } \frac{n!}{e^{\lambda }(1+\mu )^n}\frac{\lambda ^{k_1}\mu ^{(n-r)}}{k_1!r!(n-r)!}
\;=\;
\mychoose{n}{r}\, p^r (1-p)^{(n-r)}
\]
(where we use the notation $p:=1/(1+\mu )$, so that $\mu =\frac{1-p}{p}$), which
shows that this conditional marginal distribution is a binomial random
variable with parameters $n$ and
\[
p=\frac{\kappa _2\kappa _4}{\kappa _2\kappa _4+\kappa _1\kappa _2} \,.
\]

\bibliographystyle{plain}

\edo